\documentclass[11pt]{article}
\usepackage{amsmath, amsfonts, amsthm, amssymb}
\usepackage{graphicx}
\usepackage{float}
\usepackage{verbatim}
\usepackage[usenames]{color}

\numberwithin{equation}{section}
\allowdisplaybreaks

\newtheorem{thm}{Theorem}[section]
\newtheorem{lma}[thm]{Lemma}
\newtheorem{cor}[thm]{Corollary}
\newtheorem{defn}[thm]{Definition}

\newtheorem{prop}[thm]{Proposition}

\newtheorem{rem}[thm]{Remark}
\newtheorem{ques}[thm]{Question}

\renewcommand{\geq}{\geqslant}
\renewcommand{\leq}{\leqslant}
\renewcommand{\H}{\text{H}}
\renewcommand{\P}{\text{P}}

\title{A note on the 1-prevalence of continuous images with full Hausdorff dimension}

\author{Jonathan M. Fraser\footnote{Correspondence author.}\\
\emph{Mathematics Institute, Zeeman Building,}\\ \emph{University of Warwick, Coventry, CV4 7AL, UK.}\\ \emph{e-mail: jon.fraser32@gmail.com}\\ \\
James T. Hyde\\
\emph{Mathematical Institute, University of St Andrews, North Haugh,}\\
\emph{St Andrews, Fife, KY16 9SS, UK.} \\ \emph{e-mail: jth4@st-andrews.ac.uk}}
\begin{document}
\maketitle
\begin{abstract}
We consider the Banach space consisting of real-valued continuous functions on an arbitrary compact metric space.  It is known that for a prevalent (in the sense of Hunt, Sauer and Yorke) set of functions the Hausdorff dimension of the image is as large as possible, namely 1.  We extend this result by showing that `prevalent' can be replaced by `1-prevalent', i.e. it is possible to \emph{witness} this prevalence using a measure supported on a one dimensional subspace.  Such one dimensional measures are called \emph{probes} and their existence indicates that the structure and nature of the prevalence is simpler than if a more complicated `infinite dimensional' witnessing measure has to be used.
\end{abstract}

\section{Introduction}

Let $X$ be a compact metric space and let $C(X)$ denote the set of real-valued continuous functions on $X$, which is a Banach space when equipped with the supremum norm, $\| \cdot \|_\infty$.  We investigate the Hausdorff and packing dimensions of three objects related to $f \in C(X)$, namely, the image
\[
f(X) \subseteq \mathbb{R},
\]
endowed with the Euclidean metric; the product space
\[
X \times f(X) \subseteq X \times \mathbb{R}, 
\]
endowed with any of the natural and equivalent product metrics; and the graph
\[
G_f = \big\{ (x,f(x)) : x \in X \big\} \subseteq X \times f(X).
\]
Rather than compute the dimensions of these objects for specific $f \in C(X)$, instead we look to find the `generic answer'.  Namely, what are the dimensions of the above sets for a generic $f \in C(X)$?  In order to do this we need a suitable notion of genericity in Banach spaces, which we obtain using the theory of \emph{prevalence}, which we discuss in Section \ref{prevalencedef}.  We denote the Hausdorff and packing dimensions by $\dim_\H$ and $\dim_\P$ and the $s$-dimensional Hausdorff and packing measures by $\mathcal{H}^s$ and $\mathcal{P}^s$.
\\ \\
Prevalence has been used extensively in the literature to study dimensional properties of generic continuous functions.  In particular, over the past 20 years there has been considerable interest in studying the prevalent dimensions of graphs of continuous functions, using various notions of dimension, see \cite{mcclure, shaw, lowerprevalent, me_horizon,  me_prevalence, BH, me_images}.  A recent paper of Bayart and Heurteaux \cite{BH}, combined with a minor observation in \cite{me_images}, provides the most general result to date.

\begin{thm}[Bayart and Heurteaux]\label{BH1} 
Suppose $X$ is a compact subset of $\mathbb{R}^d$.  If $X$ is uncountable, then the set
\[
\{ f \in C(X) : \dim_\text{\emph{H}} G_f  = \dim_\text{\emph{H}} X +1 \}
\]
is a prevalent subset of $C(X)$.  If $X$ is finite or countable, then $\dim_\text{\emph{H}} G_f = 0$ for all $f \in C(X)$.
\end{thm}

We do not need to consider the continuum hypothesis in the above theorem because a compact subset of $\mathbb{R}^d$ is either finite, countable or has cardinality continuum (see \cite{cies}, Corollary 6.2.5).  In \cite{BH} fractional Brownian motion on $X$ provided the measure which `witnessed' the prevalence.  This meant that the additional assumption that $X$ has positive Hausdorff dimension was needed to guarantee the existence of an appropriate measure to use in the energy estimates.  Interestingly, this left the case where $\dim_\H X = 0$ to consider.  Clearly if $X$ is finite or countable then the dimension of the graph is necessarily 0, so the only open case is when $X$ has cardinality continuum but is zero dimensional.  However, by considering the prevalent Hausdorff dimension of the \emph{image}, this gap was filled by \cite{me_images} where it was observed that an older result of Dougherty see \cite[Theorem 11]{doug} does the job.
\begin{thm}[Dougherty] \label{doug1}
Let $K$ be a topological Cantor space and $C_n(K)$ be the Banach space of $\mathbb{R}^n$ valued continuous functions on $K$.  Then
\[
\{ f \in C_n(K) : \text{\emph{int}}\left(f(K) \right) \neq \emptyset \}
\]
is a prevalent subset of $C_n(K)$, where $\text{\emph{int}}\left(f(K) \right)$ denotes the interior of $f(K)$.
\end{thm}

This result has the following immediate corollary, as observed in \cite{me_images}.

\begin{cor}[Dougherty, Balka-Farkas-Fraser-Hyde]\label{doug2}
Suppose $X$ is uncountable. The following properties are prevalent in $C(X)$:
\begin{itemize}
\item[(1)] $\text{\emph{int}}\left(f(X) \right) \neq \emptyset$.
\item[(2)] $\dim_\text{\emph{H}} f(X) = \dim_\text{\emph{P}} f(X) = 1$.
\item[(3)] $0 < \mathcal{H}^1(f(X)) = \mathcal{P}^1(f(X)) < \infty$.
\item[(4)] If $\dim_\text{\emph{H}} X = 0$, then $\dim_\text{\emph{H}} G_f = 1 = \dim_\text{\emph{H}} X +1$.
\end{itemize}
\end{cor}
In particular, (4) fills the gap left by Bayart and Heurteaux.  The goal of this paper is to see to what extent Theorem \ref{BH1} and Corollary \ref{doug2} can be generalised from \emph{prevalence} to \emph{1-prevalence}.  Our main result, Theorem \ref{main}, achieves this generalisation for Corollary \ref{doug2}, parts (2) and (4), and Theorem \ref{BH1}, in the case when $\dim_\text{H} X = 0$, but we are currently unable to extend Corollary \ref{doug2}, parts (1) and (3), or Theorem \ref{BH1}, in full generality.
\\ \\
Prevalence has also been used to study dimensional properties of images of continuous functions in different settings, see \cite{prevalentimages}, where an `infinite dimensional' analogue of the projection theorems of Mattila \cite{mattila}, Kaufman \cite{kaufman} and Marstrand \cite{marstrand} was obtained.

\subsection{Prevalence and $k$-prevalence} \label{prevalencedef}

`Prevalence' provides a natural measure theoretic way of describing the \emph{generic} behaviour of  elements in a Banach space.  The theory is motivated by the fact that in finite dimensional vector spaces Lebesgue measure provides a natural tool for deciding if a property is `generic'.  Specifically, if the set of elements which do not have a certain property is a Lebesgue null set then this property is `generic'. Unfortunately, when the space in question is \emph{infinite} dimensional this approach breaks down because there is no useful analogue of Lebesgue measure.  The theory of prevalence addresses this problem by providing an analogue of `Lebesgue null' and `full Lebesgue measure' in the infinite dimensional setting.  It was first introduced in the general setting of abelian Polish groups by Christensen in the 1970s \cite{christ, christ2} and later rediscovered by Hunt, Sauer and Yorke in 1992 \cite{prevalence1} and formulated for completely metrizable topological vector spaces, but for our purposes it suffices to set up the theory for Banach spaces. Note that $(C(X), \| \cdot \|_\infty)$  is infinite dimensional and so it is natural to appeal to the theory of prevalence when considering questions of genericity in $C(X)$.

\begin{defn} \label{prevalentdef}
Let $Y$ be a Banach space.  A Borel set $F \subseteq Y$ is \emph{prevalent} if there exists a Borel measure $\mu$ on $Y$ and a compact set $K \subseteq Y$ such that $0<\mu(K) < \infty$ and
\[
\mu\big(Y \setminus (F+y)\big) = 0
\]
for all $y \in Y$. The complement of a prevalent set is called a \emph{shy} set and a non-Borel set $F \subseteq Y$ is \emph{prevalent} if it contains a prevalent Borel set.
\end{defn}

The measure $\mu$ in the above definition, which `witnesses' the prevalence of $F$ can sometimes be the push forward of Lebesgue measure concentrated on a finite dimensional subset of $X$.

\begin{defn}
A $k$-dimensional subspace $P \subseteq Y$ is called a \emph{probe} for a Borel set $F\subseteq Y$ if
\[
\mathcal{L}_P \big(Y \setminus (F+y)\big) = 0
\]
for all $y \in Y$ where $\mathcal{L}_P$ denotes the natural push forward of $k$-dimensional Lebesgue measure on $P$.  The set $F$ is $k$-\emph{prevalent} if it admits a $k$-dimensional probe.
\end{defn}

Being $k$-prevalent is clearly a stronger condition than being prevalent, and being 1-prevalent is the strongest condition of all. In some sense, the existence of a probe indicates that the structure and nature of the prevalence is simpler than if a more complicated `infinite dimensional' witnessing measure has to be used.  For more details and examples of this theory see \cite{prevalence}.

\section{Results}

To avoid trivial cases we will assume from now on that $X$ is uncountable. Our main result is the following. 

\begin{thm} \label{main}
The set of functions $f \in C(X)$ for which $\dim_\text{\emph{H}} f(X) = \dim_\text{\emph{P}} f(X) = 1$ is a 1-prevalent set.
\end{thm}

The proof of this result is deferred to Section \ref{mainproofs}.  We can use Theorem \ref{main} to glean information about the prevalent properties of the product set and the graph.

\begin{cor} \label{prodgraph}
The set of functions $f \in C(X)$ for which
\[
\dim_\text{\emph{H}} \big(X \times f(X) \big) = \dim_\text{\emph{H}} X + 1,
\]
\[
\dim_\text{\emph{P}} \big(X \times f(X) \big)  = \dim_\text{\emph{P}} X + 1,
\]
\[
\max\{1,\dim_\text{\emph{H}} X\} \leq \dim_\text{\emph{H}} G_f \leq \dim_\text{\emph{H}} X + 1
\]
and
\[
\max\{1,\dim_\text{\emph{P}} X\} \leq \dim_\text{\emph{P}} G_f \leq \dim_\text{\emph{P}} X + 1
\]
is a 1-prevalent set.
\end{cor}

Corollary \ref{prodgraph} follows immediately from Theorem \ref{main} using standard properties of the dimension theory of products, see \cite{falconer}.  A further specialisation gives the following.

\begin{cor} \label{0case}
Suppose $\dim_\text{\emph{H}} X = 0$.  Then the set of functions $f \in C(X)$ for which
\[
 \dim_\text{\emph{H}} G_f  = 1
\]
is a 1-prevalent set.  If we also have $\dim_\text{\emph{P}} X = 0$, then the set of functions $f \in C(X)$ for which
\[
 \dim_\text{\emph{H}} G_f  =   \dim_\text{\emph{P}} G_f = 1
\]
is a 1-prevalent set. 
\end{cor}
The following question is left open and is of particular interest to us.

\begin{ques}
Does Theorem \ref{BH1} remain true if \emph{prevalent} is replaced by \emph{1-prevalent}?
\end{ques}

It is interesting that in the case where $X$ is zero dimensional the above result can be obtained using a probe, but when the dimension of $X$ is strictly positive it is unclear if a probe can be used.  In particular, the proofs in \cite{BH, me_prevalence} do not use a probe. There is some evidence to suggest that the answer may be yes.  In particular, it is true if Hausdorff dimension is replaced with upper or lower box dimension and some additional properties are assumed about $X$, for example convexity.  Indeed a 1-dimensional probe was used to prove prevalence in the papers \cite{lowerprevalent, me_horizon}.  In \cite{lowerprevalent}, an \emph{intertwining condition} was introduced which is sufficient to prove that a measurable dimension function is maximised for a 1-prevalent set of graphs.  This condition is satisfied for upper and lower box dimension, but it is an interesting open problem to determine its validity for Hausdorff dimension.
\\ \\
A natural approach to proving Theorem \ref{main} is to consider the $s$-energies of measures which are push forwards of a measure on $X$ which behaves well concerning dimension.  This is reminiscent of the approach used in \cite{me_prevalence, BH}.  However, as already mentioned, problems can arise if $X$ is zero dimensional making it is less clear what this well-behaved measure should be.  The main idea in this paper to overcome this is to, instead of using push forwards of a measure on $X$ chosen by considering dimensional properties of $X$ itself, first find an injection $\phi \in C(X)$ which has an image with large dimension, and thus a sensible measure supported on it, and then consider the family of push forwards of the pull back of this measure.  Thus we find a sensible measure on $X$ to push forward which is tailored to the problem at hand, by taking the pull back of a measure we know behaves well.  We believe this approach could have applications to other problems.

\begin{rem}
Since our results on dimensions of images do not depend on the metric structure of $X$, the same results apply assuming only that $X$ is an uncountable compact separable Hausdorff space which contains a closed subspace homeomorphic to the Cantor space.  This includes any uncountable, separable, completely metrizable space.
\end{rem}

\section{Proof of Theorem \ref{main}}  \label{mainproofs}

Our key tool in estimating Hausdorff dimension will be proving finiteness of certain $s$-energies; known as the \emph{potential theoretic approach}.  Potential theoretic methods provide a powerful tool for finding lower bounds for Hausdorff dimension, which is usually more difficult than finding upper bounds.  Let $s \geq 0$ and let $\mu$ be a probability measure on $\mathbb{R}^d$.  The $s$-energy of $\mu$ is defined by
\[
I_s(\mu) =  \iint \frac{d\mu(x) \, d\mu(y)}{\lvert x-y \rvert^s}.
\]
The following theorem relates the Hausdorff dimension of a set, $F$, to the $s$-energy of probability measures supported on $F$.
\begin{prop} \label{potential}
Let $F \subset \mathbb{R}^d$ be a Borel set.  If there exists a Borel probability measure $\mu$ on $F$ with $I_s(\mu)<\infty$, then $\mathcal{H}^s(F) = \infty$ and therefore $\dim_\text{\emph{H}} F \geq s$.  Conversely, if $\mathcal{H}^s(F)>0$, then for all $0<t<s$ there exists a Borel probability measure $\mu$ on $F$ with $I_{t}(\mu) < \infty$.
\end{prop}

For a proof of this result see \cite[Theorem 4.13]{falconer}.  The following lemma is our key technical result.

\begin{lma} \label{keylem}
Let $\phi \in C(X)$ be an injection satisfying $\dim_\text{\emph{H}} \phi(X) > 0$ and let $f \in C(X)$ be arbitrary.  For almost all $\lambda \in \mathbb{R}$ we have
\[
\dim_\text{\emph{H}}(f+\lambda\phi)(X) \geq \dim_\text{\emph{H}}\phi(X).
\]
\end{lma}

\begin{proof}
Let $\phi \in C(X)$ be an injection satisfying $\dim_\text{H} \phi(X) > 0$, fix $f \in C(X)$ and $t \in \big(0,\dim_\H \phi(X) \big)$.  It follows from Proposition \ref{potential} that there exists a Borel probability measure $\mu$ supported on $\phi(X)$ with finite $t$-energy, i.e.,
\begin{equation} \label{muenergy}
I_{t} (\mu) < \infty.
\end{equation}
Since $\phi$ is a bijection between $X$ and $\phi(X)$ we can define a pull back measure $\nu = \mu \circ \phi$ on $X$.  Finally, for $\lambda \in \mathbb{R}$, we let $\nu_{\lambda} = \nu \circ (f+\lambda \phi)^{-1}$ be the push forward of $\nu$ under $f+\lambda \phi$.  Since $\nu_{\lambda}$ is supported on $(f+\lambda\phi)(X)$, it suffices to show that for almost all $\lambda \in \mathbb{R}$ the measure $\nu_{\lambda}$ has finite $t$-energy, i.e.,
\[
I_{t} (\nu_{\lambda}) < \infty.
\]
For $\lambda \in \mathbb{R}$ we have the following expression for $I_{t} (\nu_{\lambda})$.
\begin{eqnarray}
I_{t} (\nu_{\lambda}) &=& \int_ {x \in (f+\lambda\phi)(X)} \int_ {y \in (f+\lambda\phi)(X)} \frac{d\nu_{\lambda}(x)  \, d\nu_{\lambda}(y)}{\lvert x-y \rvert^{t}} \nonumber \\ \nonumber \\
&=& \int_ {u \in X} \int_ {v \in X}  \frac{d\nu (u)  \, \,d\nu (v)   }{\lvert (f+\lambda \phi)(u)-(f+\lambda \phi) (v) \rvert^{t}}  \nonumber \\ \nonumber \\ 
&=& \int_ {u \in X} \int_ {v \in X}  \frac{d\nu (u)  \, \,d\nu (v) }{\lvert \big(f(u)-f(v)\big) + \lambda \big(\phi(u)-\phi(v)\big) \rvert^{t}} \label{energyexp}
\end{eqnarray}

Let $n \in \mathbb{N}$ be arbitrary.  Observe that for fixed $a,b \in \mathbb{R}$ with $b \neq 0$ and $t>0$ we have
\begin{equation} \label{simpleintegral}
\int_{-n}^n \,\frac{d\lambda}{\lvert a+\lambda b \rvert^t} \ \leq \ \int_{-a/b-n}^{-a/b+n} \,\frac{d\lambda}{\lvert a+\lambda b \rvert^t} \  = \ \int_{-n}^{n} \,\frac{d\lambda}{\lvert \lambda b \rvert^t} \  = \ 2 \, \int_{0}^{n} \,\frac{d\lambda}{\lambda^t\lvert  b \rvert^t}.
\end{equation}
It follows from (\ref{energyexp}) and Fubini's Theorem that
\begin{eqnarray*}
\int_{-n}^n \, \Big( I_{t}(\nu_{\lambda})  \Big) \, d\lambda&=& \int_{-n}^n \,  \int_ {u \in X} \int_ {v \in X}  \frac{d\nu (u)  \, \,d\nu (v) }{\lvert \big(f(u)-f(v)\big) + \lambda \big(\phi(u)-\phi(v)\big) \rvert^{t}} \ d\lambda \\ \\
&=&   \int_ {u \in X} \int_ {v \in X}  \, \int_{-n}^n \frac{ d\lambda}{\lvert \big(f(u)-f(v)\big) + \lambda \big(\phi(u)-\phi(v)\big) \rvert^{t}} \  d\nu (u)  \, \,d\nu (v)\\ \\
&\leq&  \int_ {u \in X} \int_ {v \in X}  \, 2 \int_{0}^{n} \frac{ d\lambda}{\lambda^{t} \lvert \phi(u)-\phi(v) \rvert^{t}} \  d\nu (u)  \, \,d\nu (v) \\ \\
&\,& \quad \text{by (\ref{simpleintegral}), observing that $\phi(u)-\phi(v) \neq 0$ since $\phi$ is an injection}  \\ \\
&=&  2 \, \frac{n^{1-t}}{1-t} \, \int_ {u \in X} \int_ {v \in X}  \, \frac{ d\nu (u)  \, \,d\nu (v)}{\lvert \phi(u)-\phi(v) \rvert^{t}}   \\ \\
&\leq&  \frac{2n}{1-t} \, \int_ {u \in X} \int_ {v \in X}  \, \frac{ d (\mu \circ \phi)(u)  \, \,d(\mu \circ \phi) (v)}{\lvert \phi(u)-\phi(v) \rvert^{t}} \\ \\
&=&  \frac{2n}{1-t}\,  \int_ {p \in \phi(X)} \int_ {q \in \phi(X)}  \, \frac{ d \mu (p)  \, \,d\mu (q)}{\lvert p-q \rvert^{t}}   \\ \\
&=&  \frac{2n}{1-t} \,  I_{t}(\mu) \\ \\
&<& \infty
\end{eqnarray*}
by (\ref{muenergy}).  It follows that $\dim_\H (f+\lambda\phi)(X)  \geq t$ for almost all $\lambda  \in (-n,n)$ and hence $\dim_\H (f+\lambda\phi)(X) \geq \dim_\H \phi(X)$ for almost all $\lambda \in \mathbb{R}$ since $t$ can be chosen arbitrarily close to $\dim_\H \phi(X)$ and $n \in \mathbb{N}$ was arbitrary.  Note that we may apply Fubini's Theorem above because the measures involved are finite and the final integral is finite.
\end{proof}

\begin{lma} \label{borelset}
The set
\[
\{ f \in C(X) : \dim_\text{\emph{H}} f(X) = 1 \}
\]
is a Borel subset of $C(X)$.
\end{lma}

\begin{proof}
Let $\mathcal{K}(\mathbb{R})$ denote the set consisting of all non-empty compact subsets of $\mathbb{R}$ and equip this space with the Hausdorff metric, $d_\mathcal{H}$.  It was shown in \cite{mattilamauldin} that the function $\Delta_{\H}: (\mathcal{K}(\mathbb{R}), d_\mathcal{H}) \to \mathbb{R}$ defined by
\[
\Delta_{\H} (K) = \dim_\H K
\]
is of Baire class 2, and, in particular, Borel measurable.  Let $\Lambda : C(X) \to \mathcal{K}(\mathbb{R})$ be defined by $\Lambda(f) = f(X)$ and observe that it is continuous.  It follows that
\[
(\Delta_{\H} \circ \Lambda)^{-1}(\{1\}) = \{ f \in C(X) : \dim_\H f(K) = 1 \}
\]
is a Borel set.
\end{proof}

We can now prove Theorem \ref{main}.
\begin{proof}
Since $X$ is uncountable and compact it follows that $X$ contains a closed subset, $K$, homeomorphic to the Cantor space, see \cite[Theorem 11.11]{realanalysis}.  Let $F \subset [0,1]$ be a Cantor set with Hausdorff dimension equal to 1.  It is clear that such a set exists and for an explicit construction see, for example, the `fat Cantor set' constructed in \cite{me_prevalence}.  Since all Cantor sets are homeomorphic there exists a continuous bijection $\phi$ mapping $K$ to $F$.  We can use Tietze's Extension Theorem to extend $\phi$ to a continuous function $\Phi \in C(X)$.  We claim that
\[
P:=\{\lambda  \Phi : \lambda \in \mathbb{R} \}
\]
is a probe for $A:=\{ f \in C(X) : \dim_\H f(X) = 1 \}$.  Let $\pi_P: \mathbb{R} \to P$ be defined by $\pi_P(\lambda) = \lambda\Phi$ and let $\mathcal{L}_P=\mathcal{L}^1\circ \pi_P^{-1}$ be the push forward of 1-dimensional Lebesgue measure on $P$.  Fix $f \in C(X)$.  We have
\begin{eqnarray*}
\mathcal{L}_P\big(C(X) \setminus (A + f)\big) &=& \mathcal{L}_P \big(P \setminus (A + f)\big)\\ \\
&=& (\mathcal{L}^1\circ \pi_P^{-1})\big(\lambda  \Phi : \lambda  \Phi-f \notin A \big)\\ \\
&=& \mathcal{L}^1\big(\lambda  : \dim_\H (-f + \lambda\Phi)(X) < 1\big) \\ \\
&\leq& \mathcal{L}^1\big(\lambda  : \dim_\H (-f + \lambda\phi)(K) < \dim_\H \phi(K) \big) \\ \\
&=& 0
\end{eqnarray*}
by Lemma \ref{keylem}.  Theorem \ref{main} now follows from Lemma \ref{borelset} and the fact that $A$ admits a 1-dimensional probe.
\end{proof}

\vspace{6mm}

\begin{centering}

\textbf{Acknowledgements}

\end{centering}
Most of this work was carried out while the authors were both PhD students at the University of St Andrews.  JMF acknowledges financial support from EPSRC grant EP/J013560/1 whilst at Warwick and an EPSRC Doctoral Training Grant whilst at St Andrews.  JTH acknowledges financial support from an EPSRC Doctoral Training Grant.  The authors thank Rich\'ard Balka and \'Abel Farkas for helpful discussions, especially during the writing of the paper \cite{me_images}.  They also thank M\'arton Elekes for drawing their attention to the paper \cite{doug}.

\end{document}